\documentclass[12pt,reqno]{amsart}
\usepackage{amsthm,amsfonts,amssymb,amscd}
\newtheorem{theorem}[subsection]{Theorem}
\newtheorem{proposition}[subsection]{Proposition}

\newtheorem{lemma}[subsection]{Lemma}
\newtheorem{corollary}[subsection]{Corollary}

\theoremstyle{definition}

\newtheorem{example}[subsection]{Example}

\theoremstyle{remark}
\newtheorem{remark}[subsection]{Remark}

\newcommand{\mt}[1]{\operatorname{#1}}

\newcommand{\red}{\operatorname{red}}

\newcommand{\Diff}{\operatorname{Diff}}

\newcommand{\Sing}{\operatorname{Sing}}

\newcommand{\pal}{\text{---}}

\newcommand{\SL}{\operatorname{SL}_2(\CC)}
\newcommand{\PGL}{\operatorname{PGL}_2(\CC)}

\newcommand{\aaa}{{$A^{*}$}}
\newcommand{\ab}{{$A^{**}$}}
\newcommand{\da}{{$D^{*}$}}
\newcommand{\db}{{$D^{**}$}}
\newcommand{\ea}{{$E^{*}_6$}}

\newcommand{\CC}{{\mathbb C}}

\newcommand{\ZZ}{{\mathbb Z}}

\newcommand{\PP}{{\mathbb P}}
\newcommand{\NN}{{\mathbb N}}

\newcommand{\OOO}{{\mathcal O}}

\newcommand{\cyc}[1]{\ZZ_{#1}}

\title{Mori conic bundles with a reduced log terminal boundary}
\author{Yu.~G.~Prokhorov}

\thanks{ 
Partially supported
by the Russian Foundation of Fundamental Research grant
96-01-00820  and by a grant PECO-CEI from the Ministry of Higher Education
of France}

\address{Department of Mathematics (Algebra Section), 
Moscow State  University, Moscow 117 234, Russia}

\email{prokhoro@mech.math.msu.su\qquad \\
prokhoro@math.pvt.msu.su\\}

\date{}

\begin{document}
\begin{abstract}
We study the local structure of Mori 
contractions $f\colon X\to Z$ of relative dimension one under 
an additional assumption that there exists a reduced
divisor $S$ such that $K_X+S$ is plt and anti-ample.
\end{abstract}
\maketitle

\section*{Introduction}
Let $f\colon X\to Z$ be a  Mori conic bundle, i.~e. 
$f\colon X\to Z$ is an equi-dimensional contraction from 
threefold $X$ with only terminal
singularities onto a normal surface $Z$ 
such that $-K_X$ is $f$-ample. 
We study the analytic situation, so we assume that 
$(Z\ni o)$ is an analytic  germ at a point and $X\supset C$ is a germ 
along a reduced curve such that $f^{-1}(o)_{\red}=C$.
Such a contractions were studied in \cite{Pro},  \cite{Pro1},  \cite{Pro2}.
They are interesting in his own sake and  for applications to the 
three-dimensional birational geometry
\cite{Isk}, \cite{Isk1}.  
\par 
In this paper we consider the case when 
there exists an integer (irreducible) Weil divisor $S$ such that
$f(S)$ is a curve $R\subset Z$, $K_X+S$ is purely log terminal and anti-ample. 
The reason to consider contractions with a reduced boundary is that, 
even we start with no reduced boundary components, these may appear  
after some appropriate  blow-up (cf. \cite{Sh}).
\par
To study the three-dimensional case we have to understand the structure 
of surface contractions to curves. Such  contractions are classified in 
Sect.~1. The classification can be obtained also by  purely graph theoretical
technique (see \cite[\S 11]{KeM}). We use
complements \cite{Sh} because this method is more useful for applications.
In Sect. 2 we apply these results to threefolds. The main result is that
the log canonical divisor $K_X+S$ is either $1$- or $2$-complementary 
(Theorem \ref{mmm}).  Also a partial classification
 of such contrarctions is given.
 \par\medskip
{\sc Acknowledgments.}
The final version of this paper was written during my stay at the 
University of Lille 1, to which I am grateful for 
its hospitality and support. 

\subsection*{Notation}
We will work over the field $\CC$ of complex numbers.
A {\it contraction} is a 
projective  morphism $f\colon X\to Z$ of normal
varieties such that $f_*\OOO_X=\OOO_Z$. A contraction is said to be 
{\it extremal} if $\rho(X/Z)=1$. Usually we assume that contraction is 
$K_X$-negative, i.~e. $-K_X$ is $f$-ample. By $\Diff_S(B)$ we denote the 
different (see \cite{Sh}). For definition and properties of complements we
refer to \cite[\S 5]{Sh} or \cite[Ch. 19]{Ut}.

\section{Two-dimensional log terminal contractions}
\label{1}
\subsection{Notation}
\label{notations}
Let $(S\supset C)$ be a  a germ of
normal surface $S$ with only log terminal singularities along a
 reduced curve $C$, and $(R\ni o)$  be  a smooth  curve germ.
Let  $f\colon (S,C)\to (R,o)$ be a $K_S$-negative contraction 
such that $f^{-1}(o)_{\mt{red}}=C$.  
Then it is easy to prove that  $p_a(C)=0$ and each
irreducible component of $C$ is isomorphic to $\PP^1$.
Everywhere if we do not specify 
 the opposite we will assume that $C$ is irreducible 
 (or,  equivalently,  $\rho(S/R)=1$, i.~e. $f$ is extremal). 
Let $S_{\min}\to S$ be the minimal resolution. Since the composition map
$f_{\min}\colon S_{\min}\to R$ is flat, the fiber $f_{\min}^{-1}(o)$ is a tree 
of rational curves. Therefore it is possible to define the dual graph of 
$f_{\min}^{-1}(o)$.  We will draw it in the following way: 
$\bullet$ denotes the proper transform of $C$, while
$\circ$ denotes the exceptional curve. We attach the 
selfintersection number to the corresponding vertex. 
It is clear that the proper transform of $C$ is the only 
 $(-1)$-curve into $f_0^{-1}(o)$, so we usually omit $-1$ over $\bullet$.

\begin{example}
\label{ex2}
Let $\PP^1\times\CC^1\to\CC^1$ be the natural projection. 
Consider the following action of
$\cyc{m}$ on $\PP^1_{x,y}\times\CC^1_{u}$ 
$$
(x,y;u)\longrightarrow (x,\varepsilon^q y;\varepsilon u),\qquad
\varepsilon =\exp 2\pi i/m, \qquad
(m,q)=1.
$$
Then the morphism $f\colon
S=(\PP^1\times\CC^1)/\cyc{m}\to\CC^1/\cyc{m}$ satisfies 
the conditions above. The surface $S$ has exactly two 
singular points which are of types 
 $\frac{1}{m}(1,q)$ and $\frac{1}{m}(1,-q)$.
Morphism $f$ is toric, so  $K_S$ is $1$-complementary.
One can check that the minimal resolution of $S$ has the dual graph
$$
\stackrel{-b_1}{\circ}\pal\cdots\pal\stackrel{-b_s}{\circ}\pal
\bullet\pal\stackrel{-a_r}{\circ}\pal\cdots\pal
\stackrel{-a_{1}}{\circ},
$$
where $(b_1,\dots,b_s)$ and $(a_r,\dots,a_1)$ are uniquely
defined by $n/q$ and $n/(n-q)$, respectively (see e.~g. \cite{Br}). 
\end{example}

\begin{proposition}[see also \cite{KeM} (11.5.12)]
\label{new1}
Let  $f\colon (S,C)\to (R,o)$ be a contraction as  in \ref{notations}, 
not necessary extremal. Assume that $S$ has only DuVal singularities. 
Then  $S$ is analytically isomorphic to a surface into
$\PP^2_{x,y,z}\times\CC^1_t$ which is given by one of the following equations.
\begin{enumerate}
\renewcommand\labelenumi{(\roman{enumi})}
\item $x^2+y^2+t^nz^2=0$,  then the central fiber is a reducible conic 
and $S$ has  only one singular point which is of type $A_{n-1}$,
\item $x^2+ty^2+tz^2=0$, then then the central fiber is a non-reduced conic 
and $S$ has exactly two singular points which are of type $A_{1}$,
\item $x^2+ty^2+t^2z^2=0$, then the central fiber is a non-reduced conic 
and $S$ has only one singular point which is of type $A_{3}$, 
\item $x^2+ty^2+t^nz^2=0$, $t\ge 3$ then the central fiber is a non-reduced conic 
and $S$ has only one singular point which is of type $D_{n+1}$.
\end{enumerate}
\end{proposition}
\begin{proof}
One can show that the linear system $|-K_S|$ is very ample and 
determines an embedding $S\subset \PP^2\times R$.
Then $S$ must be given by the equation $x^2+t^ky^2+t^nz^2=0$.
\end{proof}

\subsection{Construction}
\label{construction}
Notation and assumptions as in  \ref{notations}. 
Let $d$ be the index of $C$ on $S$, i. e. the smallest 
positive integer such is that $dC\sim 0$.
If $d=1$, then $C$ is a Cartier divisor and $S$ must be smooth along
$C$, because so is $C$. If $d>1$, then there exists the 
following commutative diagram
$$
\begin{CD}
\widehat{S}@>g>>S\\
@V\widehat{f}VV@VfVV\\
\widehat{R}@>h>>R,\\
\end{CD}
$$
where  $\widehat{S}\to S$ is a cyclic \'etale outside $\Sing(S)$ cover
of degree $d$ defined by $C$  and 
$\widehat{S}\to \widehat{R}\to R$ is the Stein factorization. Then
$\widehat{f}\colon \widehat{S}\to \widehat{R}$ is also a  $K_{\widehat{S}}$-negative contraction
but not necessary extremal.
By the construction the central fiber $\widehat{C}:=\widehat{f}^{-1}(\widehat{o})$
is reducible.  Therefore   $\widehat{S}$ is smooth outside $\Sing(\widehat{S})$.
We distinguish two cases.
\begin{enumerate}
\renewcommand\labelenumi{(\Alph{enumi})}
\item 
 $\widehat{C}$ is irreducible. Then $\widehat{S}$ is smooth and 
$\widehat{S}\simeq \PP^1\times\widehat{R}$. We have the case of Example \ref{ex2}.
\item 
 $\widehat{C}$ is reducible. Then the group $\ZZ_d$ permutes  
 components  of $\widehat{C}$
transitively. Since  $p_a(\widehat{C})=0$, this gives us that all the components of 
$\widehat{C}$ passes through one point, say $\widehat{P}$ and do not intersect
each other elsewhere.  The surface $\widehat{S}$ is smooth outside 
$\widehat{P}$. Note that in this case  $K_S+C$  is not  purely log terminal, 
because so is
$K_{\widehat{S}}+\widehat{C}$. 
\end{enumerate}

\begin{corollary}
\label{l1}
In notation of \ref{construction} $S$ has at most 
two singular points on $C$. 
\end{corollary}

\begin{proof}
In the case (B) these points are   $g(\widehat{P})$
and, possible,  the image of other $d$ points on $\widehat{C}$ 
with nontrivial stabilizer. 
\end{proof}

\subsection{Additional notation}
\label{addd}
In the case (B) we denote   $P:=g(\widehat{P})$  and 
if $S$ has  two singular points, let 
$Q$ be another singular  point. 
To distinguish exceptional divisors over $P$ and $Q$ in the corresponding 
Dynkin graph 
we reserve the notation $\circ$ for exceptional divisors over $P$
and  $\circleddash$  for exceptional divisors over $Q$.

\begin{corollary}
\label{l2}
In conditions above, if   $S$ has two singular points  on $C$, then
 $K_S+C$ is log terminal near $Q$.
\end{corollary}

\begin{lemma}
\label{max}
Notation as in \ref{notations}, \ref{construction} and \ref{addd}.
Let $S'\to S$ be a finite \'etale in codimension $1$ cover.
Then there exists the decomposition $\widehat{S}\to S'\to S$.
In particular, $S'\to S$ is cyclic and the preimage of $P$ on $S'$
consists of one point.
\end{lemma}
\begin{proof}
Let $S''$ be the normalization of $S'\times_S\widehat{S}$. 
Consider the Stein factorization $S''\to R'' \to R$. 
Then $S''\to R''$ is a flat generically $\PP^1$-bundle.
Therefore for the central fiber $C''$ one has $(-K_{S''}\cdot C'')=2$,
where $C''$ is reduced and it is the preimage of $\widehat{C}$.
On the other hand, $(-K_{S''}\cdot C'')=n(-K_{\widehat{S}}\cdot \widehat{C})=2n$, 
where $n$ is the degree of $S''\to \widehat{S}$. Whence $n=1$, $S''\simeq\widehat{S}$.
This gives us the assertion.
\end{proof}

\begin{lemma}
\label{m1}
Let  $f\colon S\to (R\ni o)$ be an extremal contraction as in 
\ref{notations} (with  irreducible  $C$).
Assume that  $K_S+C$ is log terminal.
Then
\begin{enumerate}
\renewcommand\labelenumi{{\rm (\roman{enumi})}}
\item
$f\colon S\to (R\ni o)$ is analytically isomorphic to
the contraction from Example \ref{ex2}. In particular 
$S$ has exactly two singular points on  $C$
which are of types 
 $\frac{1}{m}(1,q)$ and $\frac{1}{m}(1,-q)$,
\item
$K_S+C$ is $1$-complementary.
\end{enumerate}
\end{lemma}

\begin{proof}
In the construction  \ref{construction} we have the case (A).
Therefore $f$ is toric, so $K_S+C$ is $1$-complementary.
\end{proof}

The following result gives us the classifications of surface 
log terminal contractions of relative dimension $1$.

\begin{theorem}
\label{m2}
Let  $f\colon (S\supset C)\to (R\ni o)$ be an extremal contraction 
as in \ref{notations} (with  irreducible  $C$).
Then $K_S$ is $1$, $2$ or $3$-complementary.
Moreover, there are the following cases
\begin{enumerate}
\renewcommand\labelenumi{{\rm (\roman{enumi})}}
\item
{\sl Case} \aaa:
$K_S+C$ is log terminal, 
then $K_S+C$ is $1$-complementary and $f$ is toric (see  Example \ref{ex2}).
\item
{\sl Case} \da: 
$K_S+C$ is log canonical, but not log terminal,  then  $K_S+C$ is $2$-complementary
and $f$ is a  quotient of a conic bundle of type 
{\rm (ii)} of Proposition \ref{new1} by a cyclic 
group  $\cyc{2m}$ which permutes components of the central fiber and 
acts on $S$ free in codimension $1$. The minimal 
resolution of $S$ is
$$
\begin{array}{cccccccccccccc}
\stackrel{-2}{\circ}&&&&&&&&&&&&&\\
|&&&&&&&&&&&&&\\
\stackrel{-b}{\circ}&\pal&\stackrel{-b_1}{\circ}&\pal&\cdots&\pal&\stackrel{-b_s}{\circ}&
\pal&\bullet&\pal&\stackrel{-a_r}{\circleddash}&\pal&\cdots&\pal
\stackrel{-a_{1}}{\circleddash}\\
|&&&&&&&&&&&&&\\
\stackrel{-2}{\circ}&&&&&&&&&&&&&\\
\end{array}
$$  
where $s, r\ne 0$ (recall that $S$ can be non-singular outside $P$, 
so $r=0$ is also possible).
\item
{\sl Case} \ab:
$K_S$ is $1$-complementary, but  $K_S+C$ is not log canonical. 
The minimal  resolution of $S$ is
$$
\begin{array}{ccccccccccc}
\stackrel{-a_1}{\circ}&\pal&\stackrel{-a_2}{\circ}&\pal&\cdots&\pal&
\stackrel{-a_i}{\circ}&\pal&\cdots&\pal&\stackrel{-a_r}{\circ}\\
&    &     &    &      &    & |     &    &      &    &\\
&    &     &    &      &    &\bullet&    &      &    &\\
&    &     &    &      &    & |     &    &      &    &\\
&&\stackrel{-2}{\circleddash}&\pal&\cdots&\pal&
\stackrel{-2}{\circleddash}&&&&\\
\end{array}
$$
where $r\ge 4$, $i\ne 1, r$.
\item
{\sl Case} \db:
$K_S$ is $2$-complementary, but not $1$-complementary and 
 $K_S+C$ is not log canonical.  The minimal 
resolution of $S$ is 
$$
\begin{array}{cccccccccccc}
\stackrel{-a_r}{\circ}&\pal&\cdots&\pal&
\stackrel{-a_i}{\circ}&\pal&
\cdots&\pal
\stackrel{-a_{1}}{\circ}&\pal&
\stackrel{-b}{\circ}&\pal&\stackrel{-2}{\circ}\\
&&&&|&&&&&|&&\\
&&&&\bullet&&&&&\stackrel{-2}{\circ}&&\\
&&&&|&&&&&&&\\
&&&&\stackrel{-2}{\circleddash}&\pal&
\cdots&\pal\stackrel{-2}{\circleddash}&&&&\\
\end{array}
$$ 
where $r\ge 2$, $i\ne r$.
\item
{\sl Case} \ea:
$K_S$ is $3$-complementary, but not $1$- or $2$-complementary. The minimal 
resolution of $S$ is
$$
\begin{array}{ccccccccccc}
&&&&\stackrel{-3}{\circ}&&&&&&\\
&&&&|&&&&&&\\
\stackrel{-2}{\circ}&\pal&\stackrel{-2}{\circ}&\pal&\stackrel{-b}{\circ}
&\pal&\stackrel{-2}{\circ}&\pal&\stackrel{}{\bullet}&\pal&
\underbrace{\stackrel{-3}{\circleddash}
\pal\stackrel{-2}{\circleddash}\pal
\cdots\pal\stackrel{-2}{\circleddash}}_{b-2}\\
\end{array}
$$
(it is possible that $b=2$
and $Q\in S$ is non-singular).
\end{enumerate}
\end{theorem}

\begin{remark}
\begin{enumerate}
\renewcommand\labelenumi{{\rm (\roman{enumi})}}
\item
In the case \da $K_S$ can be $1$-complementary 
\begin{enumerate}
\item[{\rm a)}]
if $P\in S$ is DuVal (see \ref{new1}  {\rm (ii)}) or
\item[{\rm b)}]
if $s=0$, $a_1=\cdots=a_r=2$, $b=r+2$.
\end{enumerate}
\item
In cases \da, \ab\; and \db\; there are additional restrictions on the 
graph of the minimal resolution. For example, in the case \ab one easily 
can check
that $(\sum_{j=1}^{i-1}a_j)-(i-1)=(\sum_{j=i+1}^{r}a_j)-(r-i)$ and
$a_i=(\text{number $\circleddash$-vertices})+2$.  
\end{enumerate}
\end{remark}

\begin{proof}
If  $K_S+C$ is  log terminal, then by Lemma~\ref{m1} we have \aaa.
If $K_S+C$ is  log canonical but not log terminal, then in Construction 
\ref{construction}
$K_{\widehat{S}}+\widehat{C}$ is also log canonical but not log terminal 
\cite[2.2]{Sh}, \cite[20.4]{Ut}.
Since $\widehat{C}$ is a Cartier divisor,  $K_{\widehat{S}}$ is canonical.
Hence $\widehat{f}$ is such as in Proposition \ref{new1}, (i).
We get the case \da. To prove that $K_S+C$ is $2$-complementary
we take on the minimal resolution $S_{\min}$ of $S$ the divisor
${D}_{\min}$ with coefficients
$$
\begin{array}{cccccccccccccc}
\stackrel{1/2}{\circ}&&&&&&&&&&&&&\stackrel{1}{\odot}\\
|&&&&&&&&&&&&&|\\
\stackrel{1}{\circ}&\pal&\stackrel{1}{\circ}&\pal&\cdots&
\pal&\stackrel{1}{\circ}&\pal&\stackrel{1}{\bullet}&\pal&
\stackrel{1}{\circleddash}
&\pal\cdots&\pal&\stackrel{1}{\circleddash}\\
|&&&&&&&&&&&&&\\
\stackrel{1/2}{\circ}&&&&&&&&&&&&&\\
\end{array}
$$  
where $\odot$ corresponds to an additional incomplete curve.
One can check that $2(K_{S_{\min}}+{D}_{\min})\sim 0$, so we have an integer Weil divisor 
$D=C+B$ on $S$ such that $K_S+C+B$ is log canonical and $2(K_S+C+B)\sim 0$. 
\par
Assume that 
$K_S$  is $1$-complementary, but  $K_S+C$ is not log terminal.
Then there exists a reduced divisor $D$ such that
$K_S+D$ is  log canonical  and linearly trivial.
By our assumption  and by \cite[9.6]{K} $C\not\subset D$.
Let  $P\in S$ is a point of index $>1$.
Then $P\in C\cap D$ and again by \cite[9.6]{K}
there are two components $D_1, D_2\subset D$ passing through  $P$.
But since  $D\cdot L=2$, where $L$ is a generic fiber of $f$,
 $D=D_1+D_2$, $P\in D_1\cap D_2$ and $P$ is the only point of index 
$>1$ on $S$. 
\par

Now we claim that 
$K_S$ is $1$-, $2$- or $3$-complementary.
Suppose that  $K_S$ is not  $1$-complementary.
Then by Lemma~\ref{m1} for some  $\alpha\le 1$ the log divisor
$K_S+\alpha C$ is log canonical, but not Kawamata log terminal.
Consider the log terminal blow-up 
$\varphi\colon(\check{S},\sum E_i+\alpha\check{C})\to (S,\alpha C)$
\cite{Sh}, where $\sum E_i$ is the reduced exceptional divisor,
$\check{C}$ is the proper transform of  $C$ and 
$K_{\check{S}}+\sum E_i+\alpha\check{C}=\varphi^*(K_S+\alpha C)$
is log terminal. Applying the  $(K_{\check{S}}+\sum E_i)$-MMP to 
$\check{S}$ we obtain  on the last step the blow-up 
$\sigma\colon\widetilde{S}\to S$ with irreducible 
 exceptional divisor $E$.
Moreover,  $\sigma^*(K_S+\alpha C)=K_{\widetilde{S}}+E+\alpha\widetilde{C}$
is log canonical,
where $\widetilde{C}$ is the proper transform of  $C$ and 
$K_{\widetilde{S}}+E$ is purely log terminal and negative over  $S$.
Since   $K_{\widetilde{S}}+E+(\alpha-\epsilon)\widetilde{C}$
is anti-ample for $0<\epsilon\ll 1$, the curve    $\widetilde{C}$ 
can be contracted  in the appropriate MMP over  $R$ and this gives 
us  $(\bar{S},\bar{C})\to R$ 
with  purely log terminal 
$K_{\bar{S}}+\bar{E}$, i.~e. by Lemma~\ref{m1} $(\bar{S},\bar{E})\to R$
 is such as in Example 
\ref{ex2}.  If    $K_{\widetilde{S}}+E$ in non-negative on 
 $\widetilde{C}$, then by  \cite[4.4]{Sh1} we can pull back 
 $1$-complements from   $\bar{S}$ 
 on  $\widetilde{S}$ and then pull down them  on $S$ \cite[5.4]{Sh}.
So below we assume that  $-(K_{\widetilde{S}}+E)$ is 
ample  over  $R$.
Then by \cite[19.6]{Ut}  complements for 
$K_E+\Diff_E(0)$ can be extended on $\widetilde{S}$.
According to  \cite[3.9]{Sh} or \cite[16.6, 19.5]{Ut}  
$\Diff_E(0)=\sum_{i=1}^3(1-1/m_i)P_i$, where
for $(m_1,m_2,m_3)$  there are the following 
possibilities:  $(2,2,m)$, $(2,3,3)$, $(2,3,4)$, $(2,3,5)$.
Further, $\bar{S}$ has exactly two singular points and these are of 
type  $\frac{1}{m}(1,q)$ and $\frac{1}{m}(1,-q)$, respectively (see 
Lemma~\ref{m1}).
Since  $\widetilde{C}$ intersects  $E$ at only one point,
this point must be singular and we have two more points with $m_i=m_j$.
We get two cases
\begin{enumerate}
\renewcommand\labelenumi{{\rm \arabic{enumi})}}
\item
$(2,2,m)$, $\widetilde{C}\cap E=\{P_3\}$, 
there is  $2$-complement.
\item
$(2,3,3)$, $\widetilde{C}\cap E=\{P_1\}$,
there is  $3$-complement.
\end{enumerate}
This proves the claim.
\par
Now  assume that  $K_S$ is $2$-complementary, but not $1$-complementary and 
 $K_S+C$ is not log canonical. Then we are in the case 1), so 
 $(\widetilde{S}\ni P_1)\simeq (\widetilde{S}\ni P_2)
 \simeq (\CC^2,0)/\cyc{2}(1,1)$ and 
 $(\widetilde{S}\ni P_3)\simeq  (\CC^2,0)/\cyc{m}(1,q)$, $(m,q)=1$.
Take the minimal resolution $S_{\min}\to \widetilde{S}$ of 
$P_1, P_2, P_3\in \widetilde{S}$. 
Over $P_1$ and $P_2$ we have only single $(-2)$-curves and 
over $P_3$ we have a chain  which must intersect the  
proper transform of $\widetilde{C}$, because $\widetilde{C}$ 
passes through $P_3$.
Since the fiber of $\widetilde{S}_{\min}\to R$ over $o$ is a 
tree of rational curves,
there are no tree of them passing through one point. Whence 
 proper transforms of $E$ and $\widetilde{C}$ on $\widetilde{S}_{\min}$
are disjoint.  Moreover, the proper transform of $E$ cannot be 
a $(-1)$-curve. Indeed, othervise contracting it we get three components of
the fiber over $o\in R$  passing through one point.
It gives us that $\widetilde{S}_{\min}$ coincides with 
the minimal resolution $S_{\min}$ of $S$.
Therefore configuration of curves on $S_{\min}$ looks like that in (iv).
We have to show only that all the curves in the down part have 
the selfintersections $-2$. Indeed, contracting $(-1)$-curves 
over $R$ we obtain $\PP^1$-bundle. Each time, we contract a $(-1)$-curve,
we have the configuration of the same type. If there is a vertex with 
selfintersection $<-2$, then on some step we get the configuration
$$
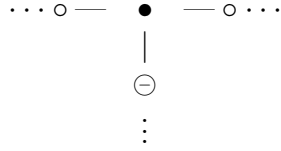
\begin{array}{ccccc}
\cdots\circ\pal&\bullet&\pal\circ\cdots\\
&|&\\
&{\circleddash}&\\
&\vdots&\\
\end{array}
$$
It is easy to see that this configuration cannot
be contracted to a non-singular point over $o\in R$, 
because contraction of the central $(-1)$-curve gives us 
configuration curves which is not a tree.
\par
The case \ea\; is very similar to \db. We omit it. 
\end{proof}

\begin{remark}
\label{2-1}
In  cases \aaa, \ab, \da\; and \db\; the canonical divisor 
$K_S$ can be $1$- or $2$-complementary 
with integer
$K_S+B$. Indeed, it is sufficient to check only in the case \db.
As in the case \da\; we take   on the minimal resolution $S_{\min}$ 
of $S$ the divisor
${B}_{\min}$ with coefficients
$$
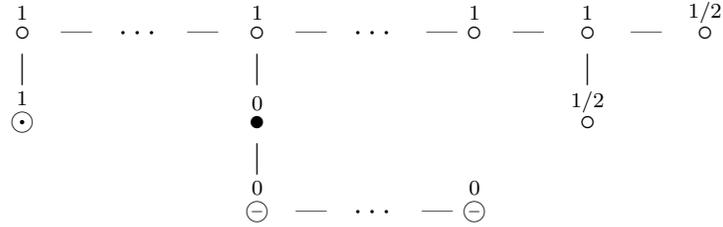
\begin{array}{cccccccccccc}
\stackrel{1}{\circ}&\pal&\cdots&\pal&
\stackrel{1}{\circ}&\pal&
\cdots&\pal
\stackrel{1}{\circ}&\pal&
\stackrel{1}{\circ}&\pal&\stackrel{1/2}{\circ}\\
|&&&&|&&&&&|&&\\
\stackrel{1}{\odot}&&&&\stackrel{0}{\bullet}&&&&&\stackrel{1/2}{\circ}&&\\
&&&&|&&&&&&&\\
&&&&\stackrel{0}{\circleddash}&\pal&
\cdots&\pal \stackrel{0}{\circleddash}&&&&\\
\end{array}
$$  
One can check that $2(K_{S_{\min}}+{B}_{\min})\sim 0$, 
so we have an integer Weil divisor 
$B$ on $S$ such that $B\cap C=\{P\}$, $K_S+B$ is log 
canonical and $2(K_S+B)\sim 0$. 
\end{remark}

\begin{corollary}
\label{dub}
In cases \da\; and \db\; there exist a double \'etale outside $P$ cover 
$g\colon S_1\to S$,  a contraction onto 
a curve $f_1\colon S_1\to R_1$ and 
 a double cover $h\colon R_1\to R$ 
such that  $hf_1=fg$. The central fiber $C_1$ of $f_1$ has exactly 
two components  which are intersects each other at $P_1:=g^{-1}(P)$.
\end{corollary}
\begin{proof}
Let $K_S+B$ be a $2$-complement from \ref{2-1}
(with integer $B$). Then $2(K_S+B)\sim 0$ and $K_S+B\not\sim 0$.
It gives us the desied cover $g\colon S_1\to S$.
The decomposition $S_1\to R_1\to R$ is the Stein factorization.
\end{proof}

\section{Three-dimensional contractions}
In this section we apply the results of the previous section to
study three-dimensional contractions with a reduced boundary. 
First we simplify the assertion  \cite[19.6]{Ut}.

\begin{lemma}
\label{prodolj}
Let $(X,S)$ be a purely log terminal pair with reduced 
$S\ne 0$  and let 
$f\colon X\to Z$ be a 
projective morphism such that $-(K_X+S)$ is nef and big. 
Given an $n$-complement $B_S^+$ of $K_S+\Diff_S(0)$, then in 
  a neighborhood of any fiber of $f$ meeting $S$ there exists an
$n$-complement $S+B^+$ of $K_X+S$ such that $B_S^+=\Diff_S(B^+)$.
\end{lemma}
\begin{proof}
By \cite[3.9]{Sh} the multiplicities $b_i$ of $\Diff_S(0)$ are standard, 
i.~e. $b_i=1-1/m_i$, $m_i\in\NN$. Consider a good resolution of 
singularities $g\colon Y\to X$. We can pull back the complements
from $S$ on its proper transform $S_Y$ by \cite[4.4]{Sh1}. 
Then applying \cite[19.6]{Ut} we get $n$-complement on $Y$.
Finally, we push down it on $X$ by \cite[5.4]{Sh}.
\end{proof}

As an easy consequence of Theorem \ref{m2} and the lemma above
 we have the following

\begin{proposition}
Let $(X\supset C)$ be a germ of a normal 
threefold  along an irreducible reduced curve, 
let  $(Z\ni o)$ be a normal surface germ and let 
$f\colon X\to Z$ be an extremal
$K_X$-negative contraction (so  $f^{-1}(o)$
is irreducible).
Assume that there exists an irreducible surface $S\subset X$ containing 
$f^{-1}(o)$
such that $f(S)$ is a curve on $Z$ and $K_X+S$ is purely log terminal
and assume that the intersection of the singular locus of
$X$ with $S$ is zero-dimensional. 
Then $K_X+S$ is $1$-,  $2$- or $3$-complementary. 
\end{proposition}

\begin{remark}
The same arguments as in the proof of Theorem~\ref{mmm}
shows that the case \ea\; is impossible, so $K_X+S$ is 
$1$- or $2$-complementary.
\end{remark}

For terminal singularities  we have more stronger results. 
In the rest of this paper we fix the following notation.

\subsection{Notation}
\label{notations2}
Let $f\colon X\to Z$ be a  Mori conic bundle, i.~e. 
$f\colon X\to Z$ is an equi-dimensional contraction from 
threefold $X$ with only terminal
singularities onto a normal surface $Z$ 
such that $-K_X$ is $f$-ample and 
$\rho(X/Z)=1$. We study the analytic situation, so we assume that 
$(S\ni o)$ is an analytic  germ at a point and $X\supset C$ is a germ 
along a reduced curve such that $f^{-1}(o)_{\red}=C$.
Since $\rho(X/Z)=1$, $C\simeq\PP^1$ \cite[1.1.1]{Pro}.
We assume that there exists an (irreducible) curve $R\subset Z$
such that $K_X+S$ is purely log terminal for $S:=f^{-1}(R)$. 
By the adjunction $S$ is normal
and has only log terminal singularities. The following fact gives us that 
 $f\colon S\to R$ is such as in 
Sect. \ref{1}.

\begin{lemma}
Let $f\colon X\to Z$, $S$ and $R$ be as in \ref{notations2}. 
Then $K_Z+R$ is purely log terminal. In particular, $R$ is non-singular.
\end{lemma}
\begin{proof}
Let $H$ be a generic hyperplane section of $X$.
By Bertini theorem $K_H+\Diff_H(S)$ is purely log terminal.
Consider the restriction $f\colon H\to Z$.
Then $\Diff_H(S)=H\cap S=f^*R$ and by the ramification formula
$K_H+\Diff_H(S)-\Delta=f^*(K_Z+R)$, where $\Delta$ is the 
ramification divisor. By \cite[\S 2]{Sh} $K_Z+R$ is purely log terminal, because so is
$K_H+\Diff_H(S)-\Delta$.
\end{proof}

Assuming $S$ is Cartier, we get a semistable contraction
(cf. e.~g. \cite{K}, \cite{Sh3})

\begin{theorem}
\label{C}
Notation as in  \ref{notations2}. 
Assume additionally that $S$ is Cartier and $f\colon X\to Z$
is not a (usual) conic bundle.
Then
\begin{enumerate}
\renewcommand\labelenumi{{\rm (\roman{enumi})}}
\item
$S$ is of type \aaa\; or \ab.
\item
$K_X+S$ is $1$-complementary.
Moreover, $K_X$ is $1$-complementary with canonical $K_X+D$.
\item
The base surface $Z$ is non-singular.
\item
$X$ contains exactly one point  of index  $>1$ and, possibly, one more 
Gorenstein terminal point.
\item
If $S$ is of type \aaa, then $f\colon S\to R$ is
 isomorphic the contraction from
Example \ref{ex2} with $n=4$. 
In particular $X$ is of index $2$ in this case 
(see Example \ref{exx} below). 
\end{enumerate}
\end{theorem}
Note that Mori conic bundles of index $2$ are completely 
classified in \cite[\S 3]{Pro1}.

\begin{proof}
By Lemma~\ref{Cartier} below $S$ has only cyclic quotient
or DuVal singulariries and by our assumption $S$
has at least one non-DuVal point. So it
is of type  \aaa, \ab\; or \da (with $s=0$, $a_1=\cdots=a_r=2$, $b=r+2$).
But the last case is also impossibe because $(P\in S)$ is not of class $T$
(see Lemma~\ref{Cartier1}). Therefore $S$ is of 
type \aaa\; or \ab.
Then by  Lemma~\ref{prodolj} there exists a $1$-complement $K_X+S+B$.
Since $S\sim 0$, $K_X+B$ is linearly trivial and canonical
(cf. \cite{Sh3}). Since $K_Z+R$ is purely log terminal and $R$ is Cartier,
$Z$ is non-singular \cite{K}. Finally, in the case  \aaa\; $S$ contains 
exactly two singular points  which are of types 
$\frac{1}{n}(1,q)$ and $\frac{1}{n}(1,n-q)$.
By Corollary \ref{rem} these singularities are of class $T$ 
iff $n=4$.
\end{proof}

\begin{lemma}[\cite{K}, \cite{KSB}]
\label{Cartier}
Let $(X\ni P)$ be a germ of a three-dimensional terminal singularity 
of index $m$ and
let $S\subset X$ be an effective Cartier 
divisor such that $K_X+S$ is purely log terminal.
Then one of the following holds
\begin{enumerate}
\renewcommand\labelenumi{{\rm (\roman{enumi})}}
\item
$(S\ni P)$ is DuVal or smooth, and then $(X\ni P)$ is cDV or  smooth,
respectively,
\item
$(S\ni P)$ is isomorphic to $\frac{1}{m^2s}(1,msm'-1)$, where $(m,m')=1$ and 
$m$ is the index of $(X\ni P)$.
\end{enumerate}
\end{lemma}

The singularities $\frac{1}{m^2s}(1,msm'-1)$, $(m,m')=1$
in \cite{KSB} are called by 
{\it singularities of class} $T$.  It is easy to see that the index of such a
singularity is equal to $m$.

\begin{corollary}
\label{rem}
If the singularity $\frac{1}{n}(1,q)$ is of class $T$,
then $n$ divides $(q+1)^2$.
\end{corollary}

Minimal resolutions of  singularities of class $T$ has very well 
description.

\begin{lemma}[\cite{KSB}]
\label{Cartier1}
\begin{enumerate}
\renewcommand\labelenumi{{\rm (\roman{enumi})}}
\item
The singularities $$
\stackrel{-4}{\circ}\qquad \text{and}\qquad
\stackrel{-3}{\circ}\pal\stackrel{-2}{\circ}\pal\stackrel{-2}{\circ}\pal
\cdots\pal\stackrel{-2}{\circ}\pal\stackrel{-3}{\circ}$$
 are of class $T$.
\item
If the singularity $\stackrel{-a_1}{\circ}\pal\cdots\pal\stackrel{-a_k}{\circ}$
is of class $T$, then so are 
$$\stackrel{-2}{\circ}\pal\stackrel{-a_1}{\circ}
\cdots\stackrel{-a_{k-1}}{\circ}\pal\stackrel{-a_k-1}{\circ}
\quad\text{and}\quad
\stackrel{-a_1-1}{\circ}\pal\stackrel{-a_2}{\circ}
\cdots\stackrel{-a_{k}}{\circ}\pal\stackrel{-2}{\circ}.$$
\item
Every singularity of class $T$ that is not DuVal can be obtained 
by starting with one of the singularities described in (i) and 
iterating the steps  described in (ii).
\end{enumerate}
\end{lemma}

Below we give an example of a semistable Mori conic bundle.

\begin{example}[cf. \cite{Pro}, 4.4]
\label{exx}
Let  $Y\subset\PP^3_{x,y,z,t}\times\CC^2_{u,v}$
is given by the equations
$$
\left\{
\begin{array}{l}
xy-z^2=ut^2\\
x^2=uy^2+v(z^2+t^2).
\end{array}
\right.
$$
Define the action of  $\cyc{2}$ on $Y$:
$$
x\to -x,\qquad y\to -y,
\qquad  z\to -z,\qquad
t\to t,\qquad  u\to u,\qquad  v\to v.
$$
Then  $f\colon X=Y/\cyc{2}\to\CC^2$
is a Mori conic bundle 
with a unique singular point of type
 $\frac{1}{2}(1,1,1)$. 
 Take a $\cyc{2}$-invariant 
 divisor $F:=Y\cap \{u=v\}$ and let  $S=F/\cyc{4}$.
Then $S\to f(S)$
 is a contraction of type \aaa. 
 More precisely, $S$ has two singular 
 points which are DuVal of type $A_3$ and a cyclic quotient of type 
 $\frac{1}{4}(1,1)$.
By the inversion of adjunction 
\cite{Sh} $K_X+S$ is  purely log terminal.
\end{example}

\begin{theorem}
\label{mmm}
Notation as in  \ref{notations2}.  Then 
\begin{enumerate}
\renewcommand\labelenumi{{\rm (\roman{enumi})}}
\item
$S$ cannot be of type \ea. 
\item
If $S$  is of type \aaa\; or \ab, then 
$K_X+S$ is $1$-complementary.
\item
If $S$  is of type \aaa, then 
$f\colon X\to Z$ is either toric (see \cite[2.1]{Pro1}) or 
is a quotient of another Mori conic bundle of index $2$
such as in  Theorem \ref{C} (v) by a cyclic group  which acts on $X$ 
free in codimension $2$.
\item
If $S$ is of type \da\; or \db, then $K_X+S$ is $2$-complementary.
In these cases $f\colon X\to Z$ is  
is a quotient of another Mori conic bundle $f'\colon X'\to Z'$
of index $1$  or $2$
and a reducible central fiber by a cyclic group of even order 
which acts on $X$ 
free in codimension $2$. Moreover
the preimage $S'$ of $S$ on $X'$ is normal and has the 
following type of the minimal resolution
\par
\begin{tabular}{ll}
Case \da.&
$\stackrel{-2}{\circ}\pal{\bullet}\pal
\stackrel{-3}{\circ}\pal\stackrel{-2}{\circ}\pal
\cdots\pal\stackrel{-2}{\circ}\pal\stackrel{-3}{\circ}
\pal{\bullet}\pal\stackrel{-2}{\circ}$,\\
Case \db.&
$
\begin{array}{ccccc}
\bullet&&&&\bullet\\
|&&&&|\\
\stackrel{-3}{\circ}&\pal&\stackrel{-2}{\circ}\pal
\cdots\pal\stackrel{-2}{\circ}&\pal&\stackrel{-3}{\circ}\\
|&&&&|\\
\bullet&&&&\bullet\\
\end{array}
$\\
\end{tabular}
\end{enumerate}
\end{theorem}

\begin{proof}
Consider the base change  \cite[3.1]{Pro}, \cite[1.9]{Pro1}
$$
\begin{CD}
X'@>g>>X\\
@V{f'}VV@V{f}VV\\
Z'@>h>>Z\\
\end{CD}
$$
where $Z'$ is smooth and $g$, $h$ are quotient morphisms by 
a finite group, say $G_0$. Let $S'$ be the preimage of $S$.
Then $S'$ is Cartier and $K_{X'}+S'$ is purely log terminal \cite[\S 2]{Sh}, 
\cite[20.4]{Ut}. By construction
the action of $G_0$  on $Z'$ is free outside $o'$ and the 
(non-singular) curve $R':=f'5S')$ is invariant under this action.
Therefore $G_0$ is cyclic.
If $X'$ has only points of index $1$, then $f'\colon X'\to Z'$ 
is a usual conic bundle (possible singular). Such quotients are described
in \cite{Pro1}. In our situation the central fiber of $f'$ cannot be multiple, 
so $f'$  either is smooth or has a reduced reducible central fiber.
From  \cite[Theorem 2.4]{Pro1} we get that  $f\colon X\to Z$ must be 
toric \cite[2.1]{Pro1} ($S$ is of type \aaa) or 
it is such as in example \ref{exd} ($S$ is of type  \da).
All the assertions hold in these cases.

Below we assume that $X'$ has at least one point of index $>1$
(and so is $S'$). 
By Lemma~\ref{max} 
the preimage of $P$ consists of one point, say $P'$.
This gives as the decomposition $(\CC^2\ni 0)\to (S'\ni P')\to (S\ni P)$
and the corresponding exact sequence
$$
1\to G_1\to G \to G_0\to 1,
$$
where $(S\ni P)\simeq (\CC^2\ni 0)/G$ and $(S'\ni P')\simeq (\CC^2\ni 0)/G_1$. 
By the construction $G_0$ is cyclic. Whence $G_1\supset [G,G]$.
On the other hand,  by Lemma~\ref{Cartier} the group
$G_1$ either is cyclic or $G_1\subset\SL$.

\subsection*{Case \ea}
In this case  $(P\in S)\simeq (\CC^2\ni 0)/G$, where the image
of $G$ in $\PGL$ is isomorphic to the alternating group 
$\mathfrak{A}_4$  \cite{Br}.
So  both $G_0$ and $G_1$ cannot be cyclic.
Therefore $(S'\ni P')$  is  a DuVal point and 
 $(X'\ni P')$ is of index $1$.
Now let $Q'_1, \dots, Q'_k$ be the preimages of $Q$. 
Since $(S\ni Q)\simeq\CC^2/\cyc{2b-3}(1,b-2)$, 
$(S'\ni Q'_i)\simeq\CC^2/\cyc{c}(1,b-2)$, where 
$c$ devides  $(2b-3)$. Thus $(1+b-2, c)=1$ and by \ref{rem}
$(S'\ni Q'_i)$ cannot be of class $T$.
Therefore $(S'\ni Q'_i)$ is smooth and $X'$ is Gorenstein, a contradiction
with our assumptions.

\subsection*{Cases \aaa\; and \ab}
Then $K_X+S$ is $1$-complementary by \ref{m2} and 
\ref{prodolj}. In the case  \aaa\;  $f'$ is (usual) conic bundle or
 such as in Theorem \ref{C} (v). This proves (ii) and (iii).

\subsection*{Cases \da\; and \db}
We claim that there is a decomposition $(X',S')\to (X_1,S_1)\to (X,S)$,
where $S_1$ is such as in Corollary~\ref{dub}. Let $K_S+B$ be an integer
$2$-complement from \ref{2-1}. It is sufficient to prove that
$K_{S'}+B'\sim 0$, where $B'$ is the preimage of $B$ on $S'$.
Assume the opposite. Note that $K_{S'}+B'$ is log canonical, but 
not purely log terminal near $P'$, since $S'\to S$ is \'etale 
in codimension $1$. Recall also that $K_{S'}+B'\sim 0$ and purely
log terminal outside $P'$. If $B'$ has locally near $P'$ two components,
then $K_{S'}+B'\sim 0$ also near $P'$, a contradiction. 
Therefore $B'$ is locally irreducible (and non-singular) near $P'$.
If $(P'\in S')$ is not DuVal, then by the classification of 
log caninical singularities with a reduced boundary \cite{K}
the minimal resolution of $(P'\in S')$ has the form
$$
\begin{array}{ccccc}
&&\stackrel{B'}{\bullet}&&\\
&&|&&\\
\stackrel{-2}{\circ}&\pal&\circ&\pal&\stackrel{-2}{\circ}\\
\end{array}
$$
But by Lemma~\ref{Cartier1} it can not be a singularity of type $T$.
Whence $(P'\in S')$ is  DuVal with the  minimal resolution of the form
$$
\begin{array}{cccccccccccccc}
\stackrel{-2}{\circ}&&&&&&&&&&&&&\\
|&&&&&&&&&&&&&\\
\stackrel{-2}{\circ}&\pal&\stackrel{-2}{\circ}&\pal&\cdots&\pal&
\stackrel{-2}{\circ}&
\pal&\stackrel{B'}{\bullet}&&&&&\\
|&&&&&&&&&&&&&\\
\stackrel{-2}{\circ}&&&&&&&&&&&&&\\
\end{array}
$$
By our assumptions $S'$ has at least one non-DuVal point.  
Then $S$ can not be of type \db, because in this case 
$S$ is DuVal at $Q'$. So $S$ is of type \da, $B'=C'$ near $P'$
and, therefore  $C'$ is irreducible. In this case it is easy to see that
$S'$ is such as in Proposition~\ref{new1} (iii)-(iv), a contradiction 
with our assumptions.
\par
Now let $(X',S')\to (X_1,S_1)\to (X,S)$ be a decomposition,
where $S_1$ is such as in Corollary~\ref{dub}.
 By our construction 
$K_{X_1}+S_1$ is $1$-complementary and so is  $K_{S_1}$.
Therefore   $S_1$ has only 
cyclic quotient  singularities.  
The group $G_0$ acts on the minimal 
 resolution of $S'$ and permutes (two) components of $C'$ in the case 
 \da\; and permutes (two) components of $B'_{S'}$ in the case \db.
 Therefore the dual graph of $S'$ is symmetric, in particular 
 the dual graph of $(P'\in S')$ is symmetric. It is easy to see that
 a singularity of class 
$T$ is symmetric iff it is as in  \ref{Cartier1} (i). But if the dual
graph of $(P'\in S')$ is $\stackrel{-4}{\circ}$ then 
$G_1$  is contained in the center of $G$. Since $G/G_1$ is 
cyclic and $G$ is not abelian, this is impossible.
Keeping in mind that in the case \da\; the fiber $C'$ has exactly two components
and in the case \db\; $K_{S'}+C'$ is not log canonical
 it is easy to see that for $S'$ we have only one of two graphs in 
 (vi).
In particular either $(P'\in S')$ and $(P'\in X')$ are of index $2$. 
This proves  our theorem.
\end{proof}

\begin{remark}
In cases  \da\; and \db\; the point 
 $(P\in S)$ is of even index. By Brieskorn's classification 
\cite{Br} $(P\in S)$ is isomorphic to $\CC^2/G$, where 
$G=(Z_{4m},Z_{2m}; D_n, C_{2n})$, $(2,m)=1$, $(n,m)=1$.
\end{remark}

Consider few examples of contractions such as in Theorem \ref{mmm}.

\begin{example}[cf. \cite{Pro},  4.1.2. (i)]
\label{exd}
Let $X'\subset\PP^2\times\CC^2$
be a hypersurface which in some coordinate system $(x,y,z;u,v)$ 
is defined by the equation
$$
x^2+y^2+(u^{2n}-v^{2n})z^2=0,
$$
Let $X'\to\CC^2$ be the natural projection and let $S'$ is the section
$X'\cap \{v=0\}$. Then the projection $S'\to R':=\{v=0\}$ 
is a contraction as in \ref{new1} (i).
Define the action of $\cyc{2}$ by
$$
u\to -u,\qquad v\to -v,\qquad x\to -x,\qquad y\to y,\qquad z\to z.
$$
$X:=X'/\cyc{2}\to\CC^2/\cyc{2}$ is a   Mori conic bundle.
$S\to R:=R'/\cyc{2}$ is a contraction of type \ab\; (if $n=1$)
or \da\;  (if $n\ge 2$). By the inversion of adjunction
 $K_X+S$ is purely log terminal. Note that 
  $S$ has a unique singular point which is DuVal so it is such as in 
  \ref{new1} (iii) or (iv). 
\end{example}

\begin{example}[cf. \cite{Pro}, 4.4]
Let  $Y\subset\PP^3_{x,y,z,t}\times\CC^2_{u,v}$
is given by the equations
$$
\left\{
\begin{array}{l}
xy=ut^2\\
z^2=u(x^2+y^2)+vt^2
\end{array}
\right.
$$
Consider the action of  $\cyc{4}$ on $Y$:
$$
x\to y,\qquad y\to -x,
\qquad  z\to iz,\qquad
t\to t,\qquad  u\to-u,\qquad  v\to-v.
$$
As above
 $f\colon Y/\cyc{4}\to\CC^2/\cyc{2}$
is a Mori conic bundle 
with a unique singular point of type
 $\frac{1}{4}(1,1,3)$. Take a $\cyc{4}$-invariant 
 divisor $F:=Y\cap \{u=v\}$ and let $S=F/\cyc{4}$. Then $S\to f(S)$
 is a contraction of type \ab. 
More precisely, $S$ has two singular 
 points which are DuVal of type $A_1$ and toric of type 
 $\frac{1}{8}(1,5)$.
\end{example}

\begin{example}[cf. \cite{Pro}, 4.3]
Let $Y\subset\PP^3_{x,y,z,t}\times\CC^2_{u,v}$ be a smooth
subvariety given by the equations
$$
\left\{
\begin{array}{l}
xy=ut^2\\
(x+y+z)z=vt^2
\end{array}
\right.
$$
Consider  the following action of $\cyc{8}$ on $Y$:
\begin{multline*}
x\to\varepsilon^{-3}z,\qquad y\to\varepsilon(x+y+z),
\qquad  z\to\varepsilon^{-3} y,\qquad
t\to t,\\
u\to \varepsilon^{-2}v,\qquad  v\to\varepsilon^{-2}u,
\end{multline*}
where $\varepsilon:=\exp(2\pi i/8)$.
We obtain a Mori conic bundle  $f\colon Y/\cyc{8}\to\CC^2/\cyc{4}$.
Let  $S:=(Y\cap \{u=v\})/\cyc{8}$. 
Then $S\to f(S)$ is a contraction of type \ab. 
More precisely, $S$ has exactly one 
 singular point which is  toric of type  $\frac{1}{16}(1,5)$.
\end{example}

Note that $K_X$ is $1$-complementary in all these examples
\cite{Pro}.

\subsection{Final remark}
Of course it is possible to consider Mori conic bundles $f\colon X\to Z$
such that  $K_X+S$-negative and 
with the condition $f(S)=Z$ instead of $\dim f(S)=1$.
But in this case  $S$ is a section of $f$ over a 
generic point, i.~e. 
$(S\cdot L)=1$ for a generic fiber $L$ of $f$.  
 By an easy construction
\cite[3.1]{Pro} or  \cite[1.9]{Pro1}
$f\colon X\to Z$ is a quotient of another 
Mori conic bundle $f'\colon X'\to Z'$ (not necessary extremal) 
over a smooth  surface $Z'$ by a cyclic group.
The preimage of $S$ is again a section  over a 
generic point. Since we are considering a germ over $o$, 
$f'$ must be smooth outside $o'$.
Then by using \cite[Lemma 4]{Isk} we can prove that $f'\colon X'\to Z'$
is smooth (see proof of Lemma 1.1 in \cite{Pro2}).
In this situation $f\colon X\to (Z\ni o)$ is analytically isomorphic to 
$\CC^2\times\PP^1/\cyc{n}(1,-1;a,0)\to\CC^2/\cyc{n}$ \cite{Pro1}.
Thus this case is not very interesting.

\end{document}